\immediate\write18{bibtex \jobname}

\documentclass[12pt, a4]{article}

\usepackage[utf8]{inputenc}
\usepackage[american]{babel}

\usepackage{geometry}
\title{Combinatorial models of global dynamics:\\ learning cycling motion from data}
\author{Ulrich Bauer, David Hien, Oliver Junge,\\  Konstantin Mischaikow, Max Snijders}
\date{January 20, 2020}

\usepackage{amsmath}
\usepackage{amssymb}
\usepackage{amsthm}
\usepackage{bm}
\usepackage{thmtools}
\usepackage{nicefrac}
\usepackage{xargs}                      
\usepackage{todonotes}
\usepackage{float}
\usepackage{graphicx}
\usepackage{comment}

\newcommand{\R}{\mathbb{R}}
\newcommand{\Z}{\mathbb{Z}}
\newcommand{\cX}{\mathcal{X}}
\newcommand{\cB}{\mathcal{B}}

\usepackage{tikz}
\usetikzlibrary{cd}
\usetikzlibrary{patterns}
\usetikzlibrary{positioning, backgrounds}
\usepackage{subcaption}
\usetikzlibrary{calc}
\usepackage[export]{adjustbox}
\usetikzlibrary{fadings}
\usetikzlibrary{decorations.markings}
\usetikzlibrary{arrows}

\newcommand{\mmbox}[1]{
	\begin{tikzpicture}
	\coordinate (cube) at (0,0);
	\drawcube{cube}{#1}
	\end{tikzpicture}	
}

\usepackage{environ}

\NewEnviron{pointofinterest}[0]{\textbf{Point of interest:} \textit{ \BODY }
}

\usepackage{hyperref}
\hypersetup{
    colorlinks=true,
    linkcolor=blue,
    filecolor=magenta,      
    urlcolor=cyan,
}
 
\definecolor{matlab1}{rgb}{0,    0.4470,    0.7410}
\definecolor{matlab2}{rgb}{0.8500,    0.3250,    0.0980}
\definecolor{matlab3}{rgb}{0.9290,    0.6940,    0.1250}
\definecolor{matlab4}{rgb}{0.4940,    0.1840,    0.5560}
\definecolor{matlab5}{rgb}{0.3010,    0.7450,    0.9330}
\definecolor{matlab6}{rgb}{0.6350,    0.0780,    0.1840}

\definecolor{macro_red}{rgb}{1, 0.34, 0}
\definecolor{macro_yellow}{rgb}{1, 0.93, 0}
\definecolor{macro_blue}{rgb}{0.17, 0, 1}
\definecolor{macro_green}{rgb}{0, 1, 0.35}
\definecolor{macro_black}{rgb}{0, 0, 0}

\definecolor{macro_greenred}{rgb}{0.5, 0.52, 0.22}
\definecolor{macro_blueyellow}{rgb}{0.49, 0.47, 0.49}
\definecolor{macro_greenyellow}{rgb}{0.5, 0.97, 0.18}
\definecolor{macro_bluered}{rgb}{0.44, 0.20, 0.49}
\definecolor{macro_bluegreen}{rgb}{0.33, 0.49, 0.52}
\definecolor{macro_redyellow}{rgb}{0.64, 0.49, 0.18}
\definecolor{macro_bluegreenredyellow}{rgb}{0.35, 0.35, 0.27}

\definecolor{julia_lightgreen}{rgb}{0.56, 0.93, 0.56}
\definecolor{julia_darkgreen}{rgb}{0.0, 0.39, 0.0}
\definecolor{julia_lightblue}{rgb}{0.68, 0.85, 0.9}
\definecolor{julia_darkblue}{rgb}{0.0, 0.0, 1.0}
\definecolor{julia_yellow}{rgb}{1.0, 1.0, 0.0}
\definecolor{julia_white}{rgb}{1.0, 1.0, 1.0}
\definecolor{julia_orange}{rgb}{1.0, 0.65, 0.0}
\definecolor{julia_purple}{rgb}{0.5, 0.0, 0.5}

\newcommand{\drawcube}[2]{
	\begin{scope}[scale=.5]
		\draw[fill=#2, line cap=round, line width=0.05mm] ($(#1)$) rectangle ($(#1) + (.5,.5)$);
		\draw[fill=#2, line cap=round, line width=0.05mm] ($(#1) + (.5, .5)$) -- ($(#1) + (.5, .5) + (.2, .1)$) -- ($(#1) + (.5, .5) + (-.3,.1)$) -- ($(#1) + (0, .5)$) -- cycle;
		\draw[fill=#2, line cap=round, line width=0.05mm] ($(#1) + (.5, .5)$) -- ($(#1) + (.5, .5) + (.2, .1)$) -- ($(#1) + (.5, .5) + (.2, .1) + (0, -.5)$) -- ($(#1) + (.5, 0)$) -- cycle;
	\end{scope}
}

\definecolor{hsv_0}{rgb}{1, 0, 0}
\definecolor{hsv_1}{rgb}{0.5000, 1.0000, 0}
\definecolor{hsv_2}{rgb}{0, 1, 1}
\definecolor{hsv_3}{rgb}{0.5, 0, 1.0}
\definecolor{hsv_4}{rgb}{1, 0, 0}

\definecolor{color_1}{rgb}{0.000000, 0.906250, 1.000000}
\definecolor{color_2}{rgb}{0.000000, 0.812500, 1.000000}
\definecolor{color_3}{rgb}{0.000000, 0.906250, 1.000000}
\definecolor{color_4}{rgb}{0.000000, 1.000000, 0.906250}
\definecolor{color_5}{rgb}{0.000000, 1.000000, 0.531250}
\definecolor{color_6}{rgb}{0.000000, 0.718750, 1.000000}
\definecolor{color_7}{rgb}{0.000000, 0.343750, 1.000000}
\definecolor{color_8}{rgb}{0.000000, 1.000000, 0.156250}
\definecolor{color_9}{rgb}{0.218750, 1.000000, 0.000000}
\definecolor{color_10}{rgb}{0.312500, 1.000000, 0.000000}
\definecolor{color_11}{rgb}{0.406250, 1.000000, 0.000000}
\definecolor{color_12}{rgb}{0.406250, 1.000000, 0.000000}
\definecolor{color_13}{rgb}{0.593750, 1.000000, 0.000000}
\definecolor{color_14}{rgb}{0.968750, 1.000000, 0.000000}
\definecolor{color_15}{rgb}{0.031250, 0.000000, 1.000000}
\definecolor{color_16}{rgb}{0.406250, 0.000000, 1.000000}
\definecolor{color_17}{rgb}{0.593750, 0.000000, 1.000000}
\definecolor{color_18}{rgb}{0.593750, 0.000000, 1.000000}
\definecolor{color_19}{rgb}{0.687500, 0.000000, 1.000000}
\definecolor{color_20}{rgb}{0.781250, 0.000000, 1.000000}
\definecolor{color_21}{rgb}{1.000000, 0.000000, 0.843750}
\definecolor{color_22}{rgb}{1.000000, 0.656250, 0.000000}
\definecolor{color_23}{rgb}{1.000000, 0.281250, 0.000000}
\definecolor{color_24}{rgb}{1.000000, 0.000000, 0.468750}
\definecolor{color_25}{rgb}{1.000000, 0.000000, 0.093750}
\definecolor{color_26}{rgb}{1.000000, 0.093750, 0.000000}
\definecolor{color_27}{rgb}{1.000000, 0.187500, 0.000000}
\definecolor{color_28}{rgb}{1.000000, 0.093750, 0.000000}

\definecolor{circle_1}{rgb}{1.000000 0.000000 0.000000}
\definecolor{circle_2}{rgb}{1.000000 0.093750 0.000000}
\definecolor{circle_3}{rgb}{1.000000 0.187500 0.000000}
\definecolor{circle_4}{rgb}{1.000000 0.281250 0.000000}
\definecolor{circle_5}{rgb}{1.000000 0.375000 0.000000}
\definecolor{circle_6}{rgb}{1.000000 0.468750 0.000000}
\definecolor{circle_7}{rgb}{1.000000 0.562500 0.000000}
\definecolor{circle_8}{rgb}{1.000000 0.656250 0.000000}
\definecolor{circle_9}{rgb}{1.000000 0.750000 0.000000}
\definecolor{circle_10}{rgb}{1.000000 0.843750 0.000000}
\definecolor{circle_11}{rgb}{1.000000 0.937500 0.000000}
\definecolor{circle_12}{rgb}{0.968750 1.000000 0.000000}
\definecolor{circle_13}{rgb}{0.875000 1.000000 0.000000}
\definecolor{circle_14}{rgb}{0.781250 1.000000 0.000000}
\definecolor{circle_15}{rgb}{0.687500 1.000000 0.000000}
\definecolor{circle_16}{rgb}{0.593750 1.000000 0.000000}
\definecolor{circle_17}{rgb}{0.500000 1.000000 0.000000}
\definecolor{circle_18}{rgb}{0.406250 1.000000 0.000000}
\definecolor{circle_19}{rgb}{0.312500 1.000000 0.000000}
\definecolor{circle_20}{rgb}{0.218750 1.000000 0.000000}
\definecolor{circle_21}{rgb}{0.125000 1.000000 0.000000}
\definecolor{circle_22}{rgb}{0.031250 1.000000 0.000000}
\definecolor{circle_23}{rgb}{0.000000 1.000000 0.062500}
\definecolor{circle_24}{rgb}{0.000000 1.000000 0.156250}
\definecolor{circle_25}{rgb}{0.000000 1.000000 0.250000}
\definecolor{circle_26}{rgb}{0.000000 1.000000 0.343750}
\definecolor{circle_27}{rgb}{0.000000 1.000000 0.437500}
\definecolor{circle_28}{rgb}{0.000000 1.000000 0.531250}
\definecolor{circle_29}{rgb}{0.000000 1.000000 0.625000}
\definecolor{circle_30}{rgb}{0.000000 1.000000 0.718750}
\definecolor{circle_31}{rgb}{0.000000 1.000000 0.812500}
\definecolor{circle_32}{rgb}{0.000000 1.000000 0.906250}
\definecolor{circle_33}{rgb}{0.000000 1.000000 1.000000}
\definecolor{circle_34}{rgb}{0.000000 0.906250 1.000000}
\definecolor{circle_35}{rgb}{0.000000 0.812500 1.000000}
\definecolor{circle_36}{rgb}{0.000000 0.718750 1.000000}
\definecolor{circle_37}{rgb}{0.000000 0.625000 1.000000}
\definecolor{circle_38}{rgb}{0.000000 0.531250 1.000000}
\definecolor{circle_39}{rgb}{0.000000 0.437500 1.000000}
\definecolor{circle_40}{rgb}{0.000000 0.343750 1.000000}
\definecolor{circle_41}{rgb}{0.000000 0.250000 1.000000}
\definecolor{circle_42}{rgb}{0.000000 0.156250 1.000000}
\definecolor{circle_43}{rgb}{0.000000 0.062500 1.000000}
\definecolor{circle_44}{rgb}{0.031250 0.000000 1.000000}
\definecolor{circle_45}{rgb}{0.125000 0.000000 1.000000}
\definecolor{circle_46}{rgb}{0.218750 0.000000 1.000000}
\definecolor{circle_47}{rgb}{0.312500 0.000000 1.000000}
\definecolor{circle_48}{rgb}{0.406250 0.000000 1.000000}
\definecolor{circle_49}{rgb}{0.500000 0.000000 1.000000}
\definecolor{circle_50}{rgb}{0.593750 0.000000 1.000000}
\definecolor{circle_51}{rgb}{0.687500 0.000000 1.000000}
\definecolor{circle_52}{rgb}{0.781250 0.000000 1.000000}
\definecolor{circle_53}{rgb}{0.875000 0.000000 1.000000}
\definecolor{circle_54}{rgb}{0.968750 0.000000 1.000000}
\definecolor{circle_55}{rgb}{1.000000 0.000000 0.937500}
\definecolor{circle_56}{rgb}{1.000000 0.000000 0.843750}
\definecolor{circle_57}{rgb}{1.000000 0.000000 0.750000}
\definecolor{circle_58}{rgb}{1.000000 0.000000 0.656250}
\definecolor{circle_59}{rgb}{1.000000 0.000000 0.562500}
\definecolor{circle_60}{rgb}{1.000000 0.000000 0.468750}
\definecolor{circle_61}{rgb}{1.000000 0.000000 0.375000}
\definecolor{circle_62}{rgb}{1.000000 0.000000 0.281250}
\definecolor{circle_63}{rgb}{1.000000 0.000000 0.187500}
\definecolor{circle_64}{rgb}{1.000000 0.000000 0.093750}

\newcommand{\cK}{\mathcal{K}}

\begin{document}

\maketitle

\begin{abstract}
We describe a computational method for constructing a coarse combinatorial model of some dynamical system in which the macroscopic states are given by elementary cycling motions of the system.  Our method is in particular applicable to time series data.  We illustrate the construction by a perturbed double well Hamiltonian as well as the Lorenz system.
\end{abstract}

\section{Introduction}

Conley's fundamental theorem \cite{Co78} characterizes the global structure of the dynamics of a continuous map on a compact metric space. It states that the space can be decomposed into a (chain) recurrent set and its complement, on which the map behaves gradient-like, i.e.\ trajectories transit from one recurrent component to another.  Around the turn of the century, a computational approach to this theory has been developed  \cite{Os99,Gaio2002,Mi02,KaMiVa05,BaKa06,MiMrWe16}.  

Relatedly, ideas have been put forward in order to characterize the dynamics \emph{within} a transitive component of the chain recurrent set.  For example, in \cite{DeJu99}, certain eigenfunctions of the transfer (or push forward) operator have been used in order to decompose a transitive component into, e.g., \emph{almost invariant} (aka metastable) subsets.  

The purpose of this note is to outline a computational procedure by which certain \emph{cycling behaviour} of the system can be detected and agglomerated into a coarse model.   More precisely, we describe how to detect whether the system exhibits motions along a topological circle in some geometric complex that represents a transitive recurrent component of the system.  

In particular, our technique is applicable if no model is available, but the dynamics is only given in form of a time series of data points $x_k=x(t_k)\in\R^d $, $k=1,\ldots,m$, that are, e.g., sampled from solution curves $x:[0,1]\to\R^d$ of some differential equation or constructed by a time-delay embedding of scalar measurement data. In cases where this data set is large, e.g., when the are multiple time scales in the system, a straightforward construction of a complex with this many points will be computationally infeasible.  We propose to preprocess the data by quantizing it, yielding a complex whose size essentially scales with the dimension of the underlying recurrent set.  

\section{The construction}
\label{sec:construction}
Given the time series $x_1,\dots,x_m $ in $\R^d$, we 
construct a combinatorial model which captures different types of cycling motion. Our pipeline consists of three main steps: preprocessing, finding dynamically relevant coordinates and construction of a combinatorial model.

In the preprocessing step we construct a quantization of the time series resp.\ the associated point cloud $X = \{ x_1,\dots,x_m\}$ by projecting onto a suitably chosen cubical grid which leads to a considerably reduction of the data to be processed.  In a second step, we use the topology of a Vietoris--Rips complex constructed on the quantized point cloud $X$ to obtain a set of coordinates which captures cycling motion of the time series. In the final step, we use the coordinates of the previous step to construct a combinatorial macro model for the dynamics.

\subsection{Quantization}

\begin{figure}[H]
	\centering
	\begin{subfigure}[t]{0.45\textwidth}
		\centering
		\includegraphics[width=7cm,height=4cm]{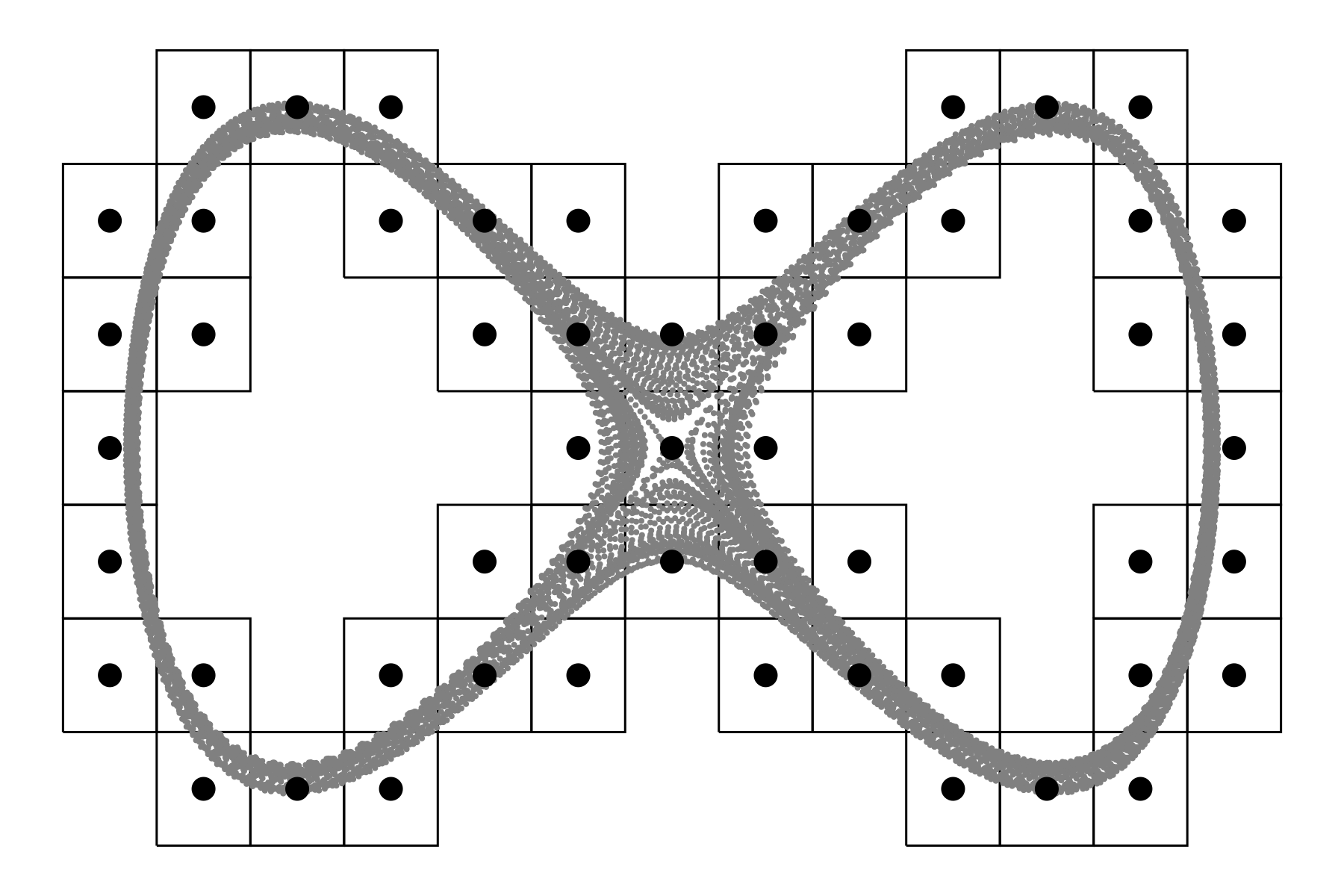}	
		\caption{time series and cubical quantization}
		\label{fig:quantization}
	\end{subfigure}\hfill
	\begin{subfigure}[t]{0.45\textwidth}
		\centering
		\includegraphics[width=7cm,height=4cm]{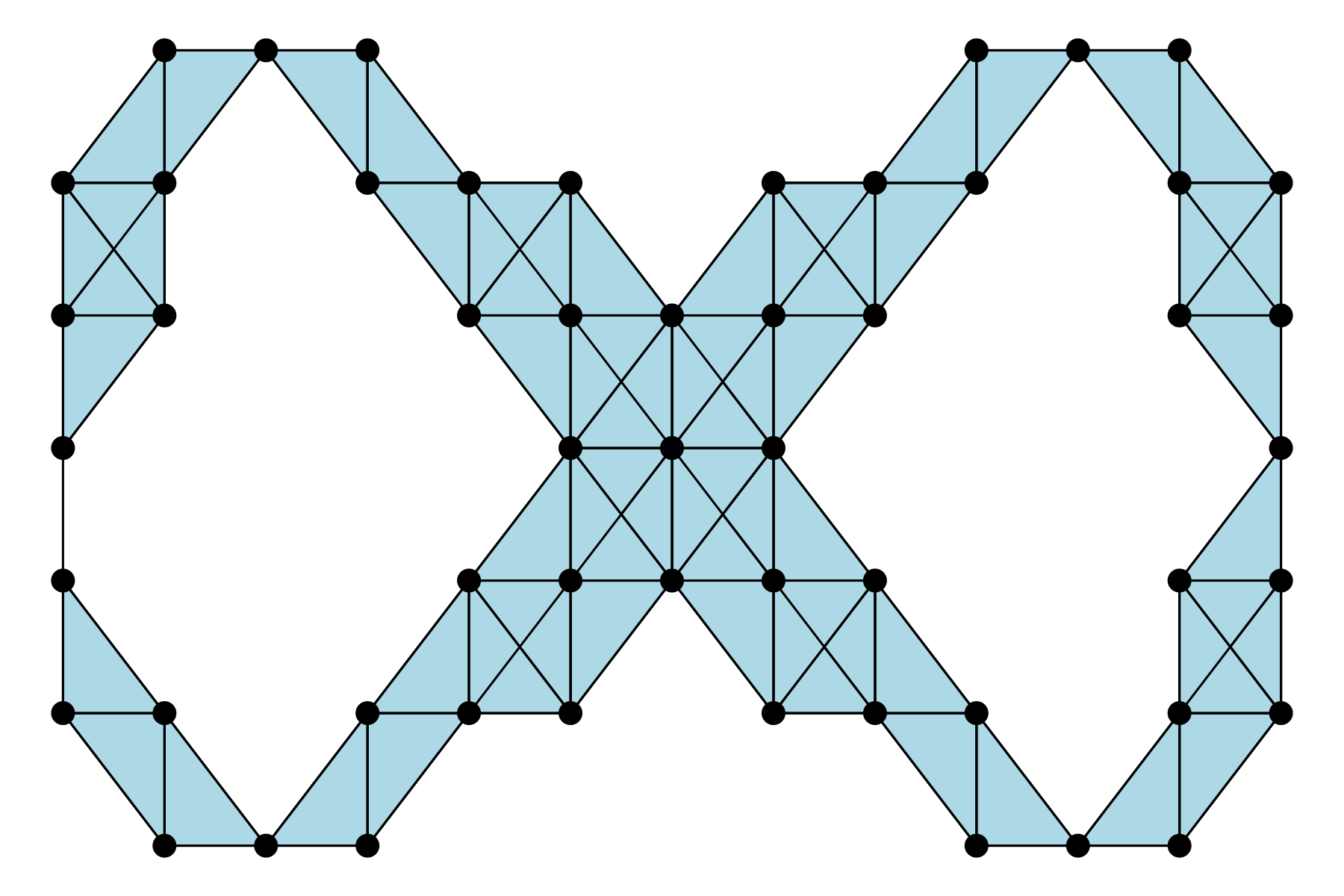}	
		\caption{Vietoris--Rips complex}
		\label{fig:Rips_complex}
	\end{subfigure}
	\caption{Time series, its quantization and the resulting Vietoris--Rips complex.}
	\label{fig:trajectory_and_covering}
\end{figure}

In the first step, the time series is quantized.  To this end, we choose a radius\ $r\in\R_{>0}^d$ and consider the \emph{cubical grid}\[
\cB = \cB(r) = \left\{ \prod_{\ell=1}^d r_\ell\left[z_\ell-\tfrac12,z_\ell+\tfrac12\right) \mid z\in \Z^d\right\}.
\]
Since the elements of $\cB$ (which we call \emph{cubes} or \emph{boxes}) form a partition of $\R^d$, we can define $Q:\R^d\to\cB$ by mapping each point to the unique cube containing the point. Then
\[
\cX := \{Q(x) \mid x\in X\}
\]
is a \emph{cubical} or \emph{box covering} of the point cloud $X$.  For a cube $\xi = \prod_{\ell=1}^d r_\ell\left[z_\ell-\tfrac12,z_\ell+\tfrac12\right)$ let $z(\xi) = (z_1,\ldots,z_d)\in\Z^d$ be its \emph{center}.  We can identify $\cX$ with the subset
\[
Z := \{z(\xi)\mid \xi \in\cX\}
\]
of the integer lattice $\Z^d$. The set $Z$ of box centers is called the \emph{quantization of the point cloud} $X$. Fig.~\ref{fig:quantization} shows a time series with its cubical cover and the corresponding set of box centers $Z$.

We then resample the time series such that consecutive points lie in different cubes. For this, we set $\tau(1) = 1$, recursively define
\[
\tau(i) = \min \{ j > \tau(i-1) \mid Q(x_j) \neq Q(x_{\tau(i-1)}) \}
\]
and set $ \hat{x}_i = z\circ Q(x_{\tau(i)}) $, $i \in [1,T]$ to be the \emph{quantization} of the time series $ x_1,\dots,x_m $. Here we let $T$ denote the largest finite value of $\tau$ and $ [1,T] := \{1,\dots,T\} $.

\subsection{Coordinates for dynamics}

After reviewing the construction of cyclic coordinates from cohomology and how these coordinates are lifted for a time series, we explain how to find elements of $ H^1(\cK,\Z) $ that induce dynamically relevant coordinates.

\paragraph{1. Cohomology and circular coordinates.}

Given a set of box centers $Z$, we construct the Vietoris--Rips complex $ \cK = \operatorname{VR}(Z,d,\delta) $, where $ d(x,y) = \| x-y \|_\infty $ and $\delta = 1$. Note that the choice of $d$ and $\delta$ allows a point in $Z$ to be connected to all its diagonal neighbors. An example is shown in Fig.~\ref{fig:Rips_complex}.

Next, we compute a basis $ B $ of $ H^1(\cK,\Z)$. Using the procedure introduced in \cite{DeSilva2011}, a circle valued coordinate $\theta_\alpha$ can be constructed for each generator $\alpha$ in $B$. More precisely, as a function on the vertices, the coordinate $ \theta_{\alpha} : Z \rightarrow S^1$  can be chosen as any solution of the optimization problem
\begin{equation*}
	\operatorname{argmin} \{ \| v + d_0\theta \|^2 \mid \theta \in C^0(X,\R) \}
\end{equation*}
composed with the canonical projection $\pi_{S^1} : \R \rightarrow S^1=\R/\Z$. Here $d_0$ denotes the coboundary operator, $\| v\|^2 $ is the sum of $ v(e)^2 $ over all edges $e$ in $\cK$ and $ v\in \alpha $ is any 1-cocycle in the respective coset. We remark that $ \theta_\alpha $ does not depend on the choice of $v$ and is unique up to an additive constant on each connected component of $\cK$.

\begin{figure}[h]
	\centering
	\begin{subfigure}[t]{0.45\textwidth}
		\centering
		\includegraphics[width=.8\textwidth]{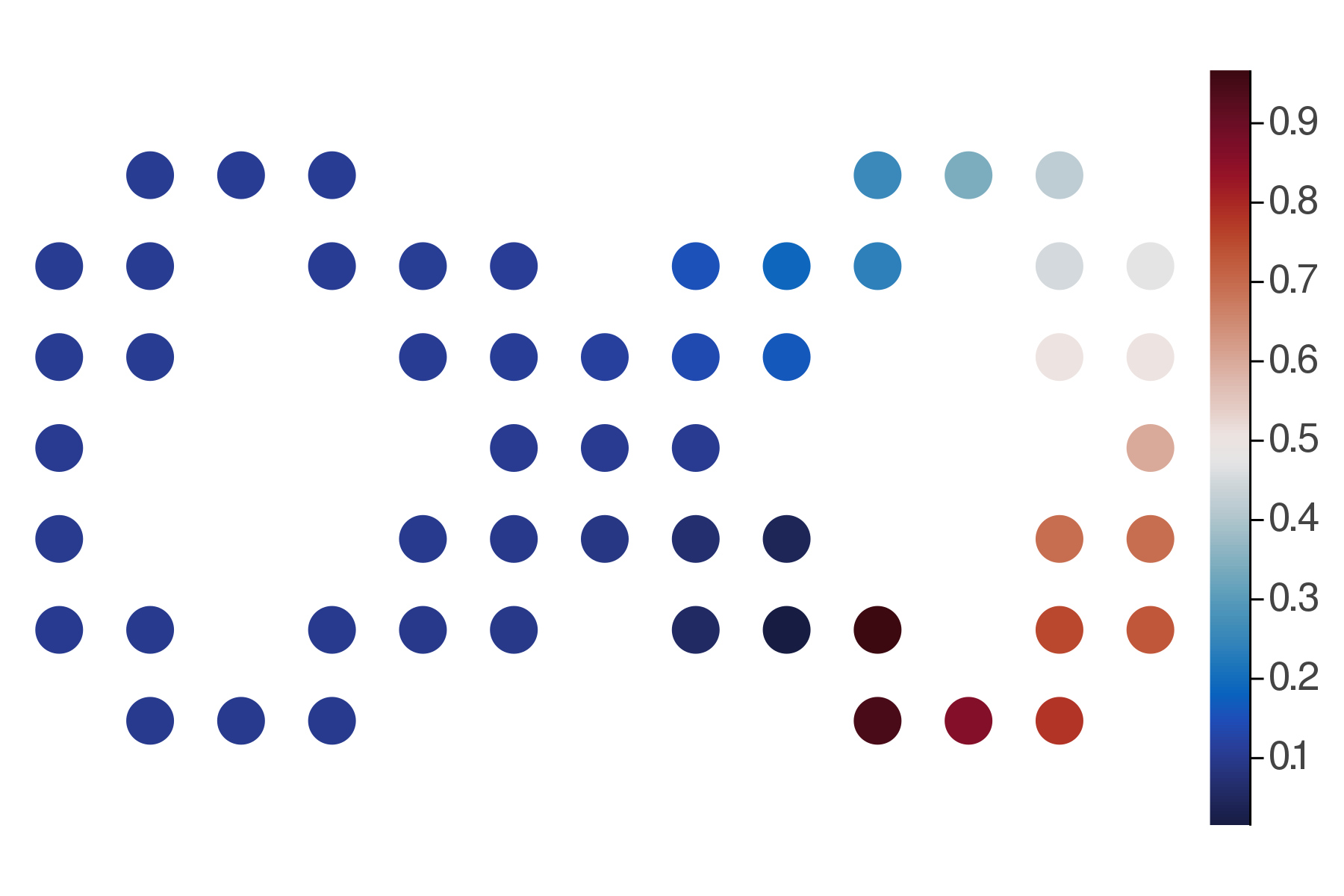}	
		\caption{Coordinate which describes the right hole}
		\label{fig:coord-good1}
	\end{subfigure}\hfill
	\begin{subfigure}[t]{0.45\textwidth}
		\centering
		\includegraphics[width=.8\textwidth]{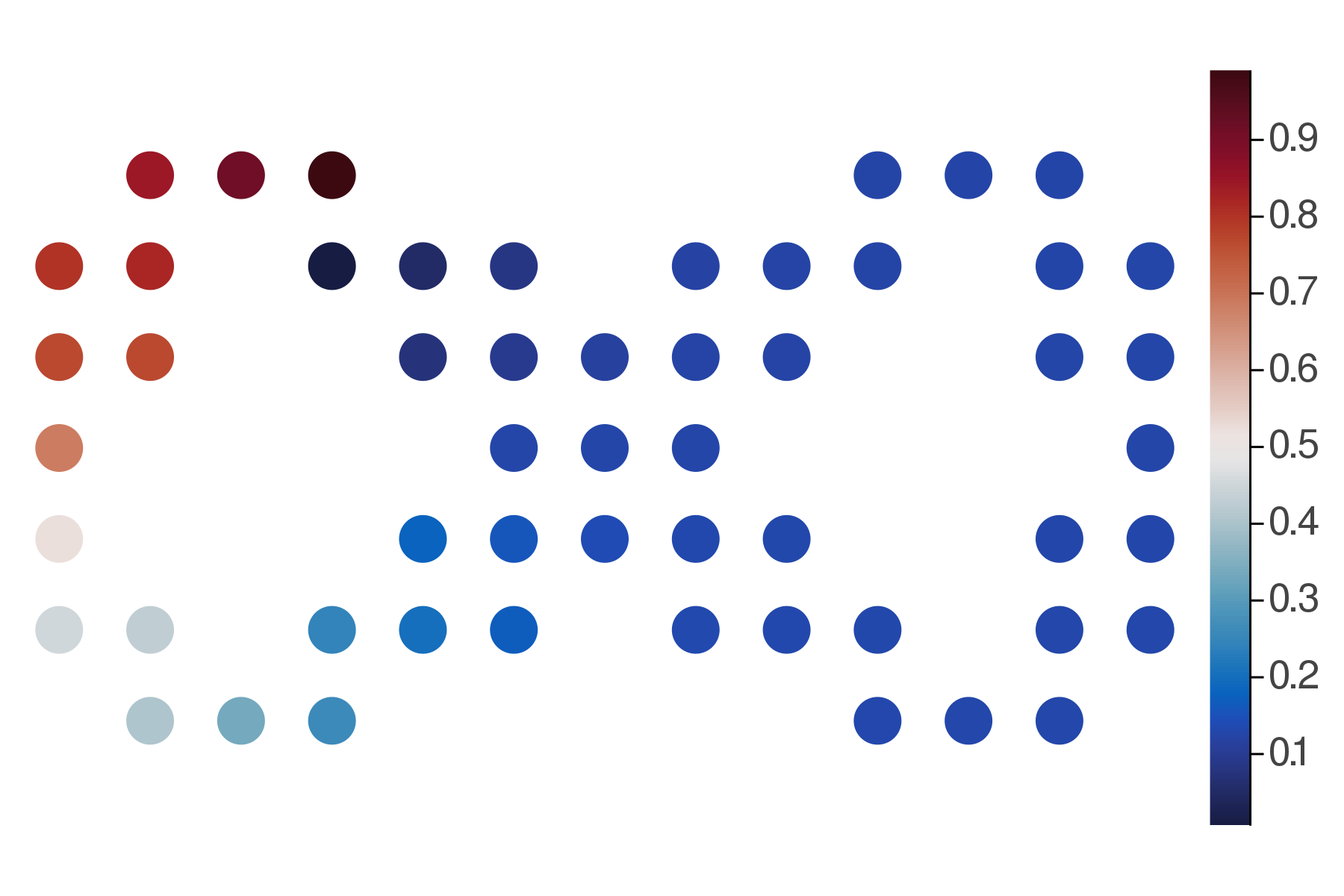}	
		\caption{Coordinate which desbribes the left hole}
		\label{fig:coord-good2}
	\end{subfigure}
	\caption{Circular coordinates for the complex in Fig.~\ref{fig:Rips_complex}. Since the complex contains two holes, its first cohomology is generated by two 1-cocycles and we compute two circular coordinates. The coordinates in this figure are particularly nice since both of them capture cycling around one of the holes.}
	\label{fig:coord-good}
\end{figure}

\begin{figure}[h]
	\centering
	\begin{subfigure}[t]{0.45\textwidth}
		\centering
		\includegraphics[width=.8\textwidth]{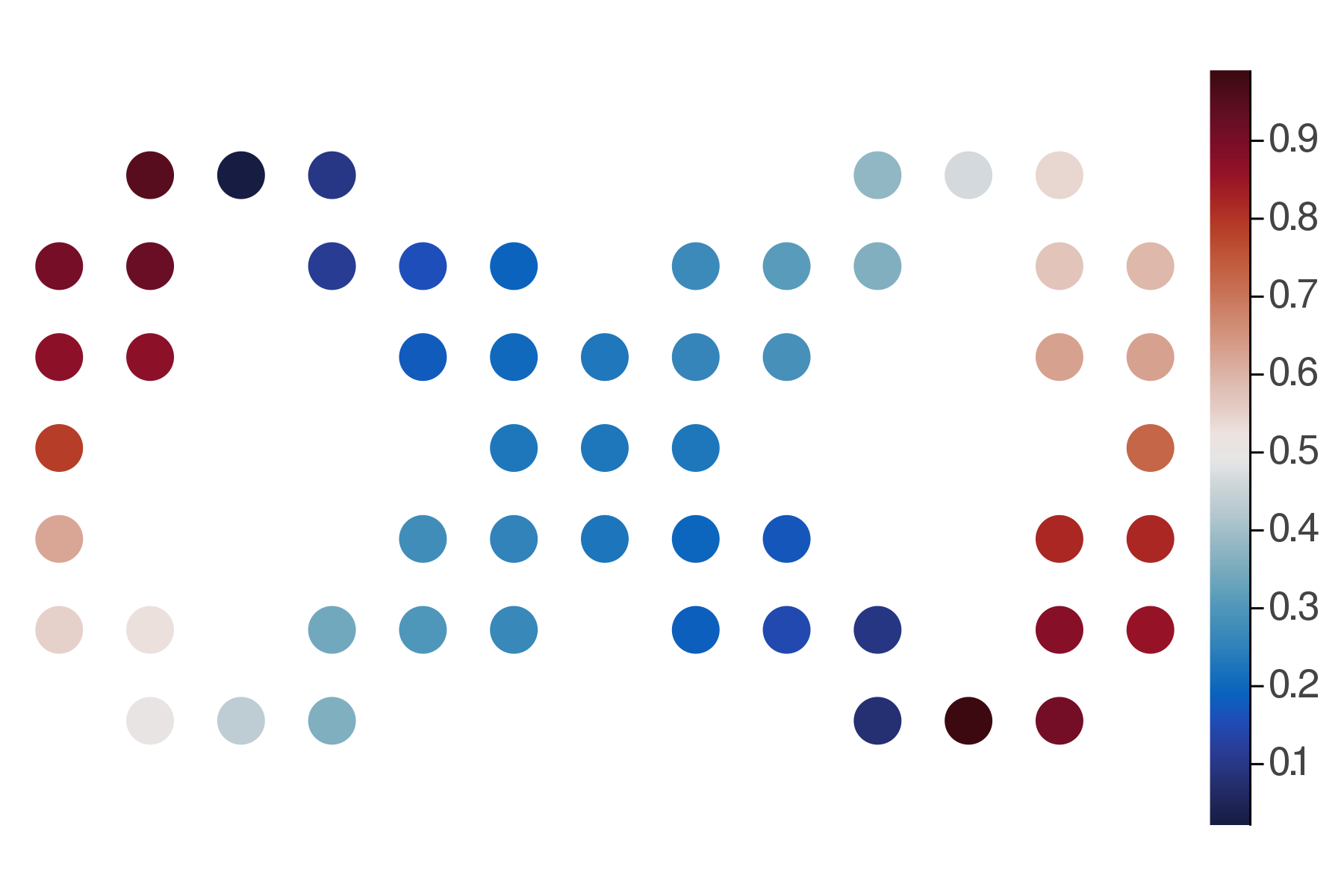}
		\caption{This coordinate varies around both holes.}
		\label{fig:coord-bad1}
	\end{subfigure}\hfill
	\begin{subfigure}[t]{0.45\textwidth}
		\centering
		\includegraphics[width=.8\textwidth]{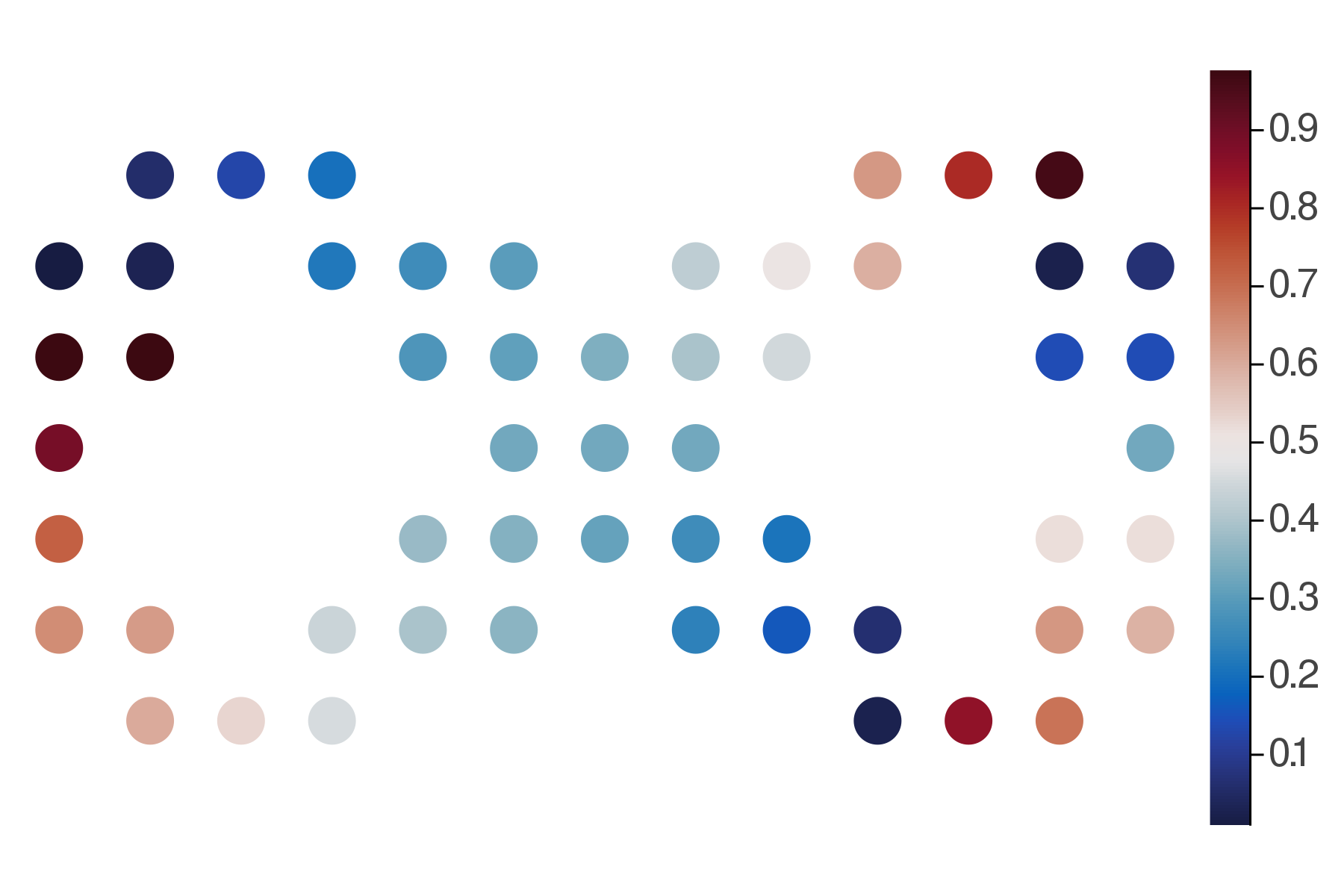}
		\caption{This coordinate varies around both holes. Additionally, it maps half-turns around the right hole surjectively onto $S^1 $.}
		\label{fig:coord-bad2}
	\end{subfigure}
	\caption{A pair of bad coordinates.}
	\label{fig:coord-bad}
\end{figure}

The set of coordinates $ \{\theta_\alpha\}_{\alpha \in B} $ however does depend on the choice of the basis $B$. As an example, consider the complex in Fig.~\ref{fig:Rips_complex}. Since its first cohomology group is isomorphic to $ \Z\times\Z $, each basis consists of two cocycles. Figures \ref{fig:coord-good} and \ref{fig:coord-bad} show coordinates for two different bases. The problem of choosing a suitable basis is addressed in a later section.

\begin{figure}[h]
	\centering
	\begin{subfigure}[t]{0.45\textwidth}
		\centering
		\includegraphics[width=.9\textwidth]{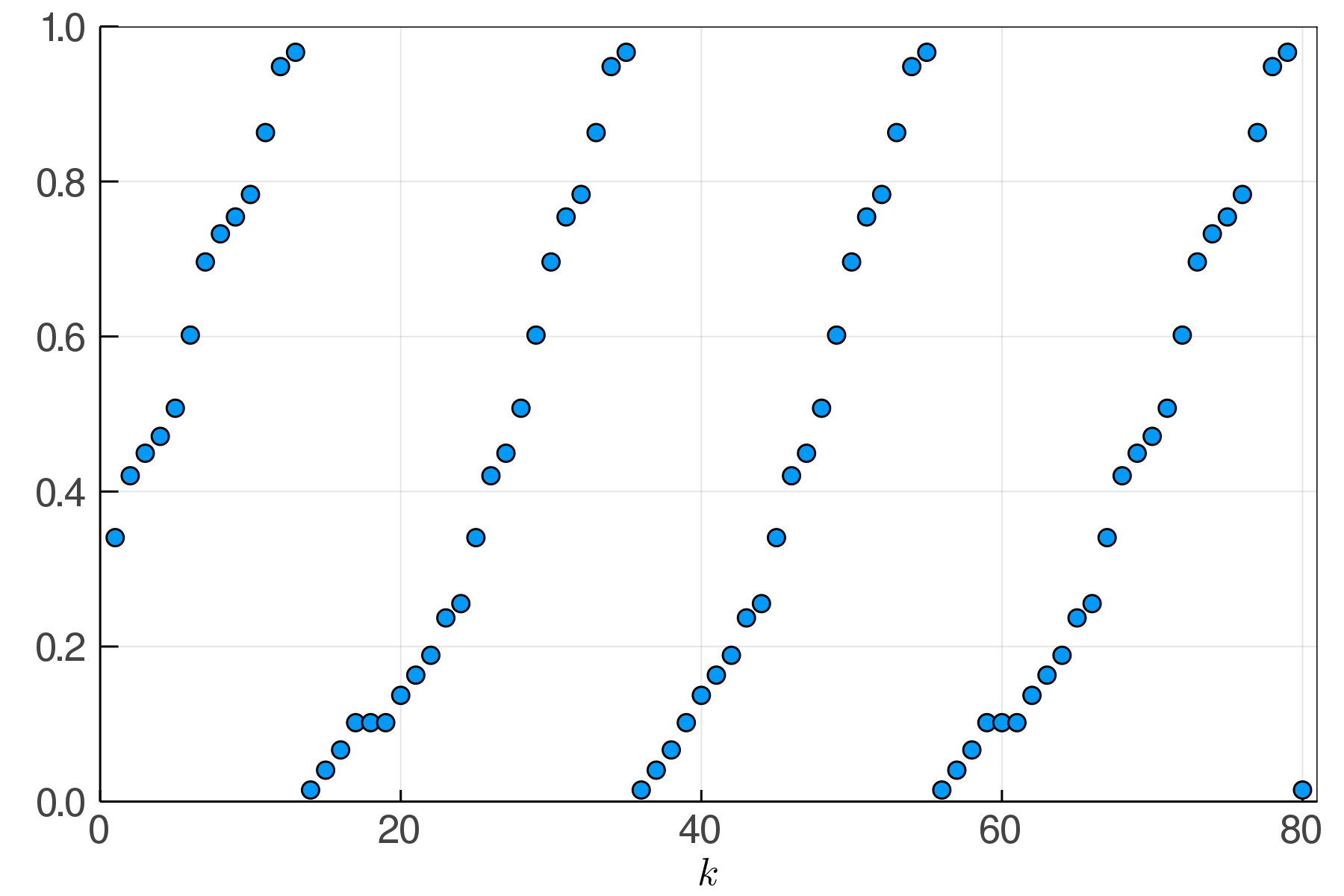}	
		\caption{Coordinate change over time.}
		\label{fig:coord-time}
	\end{subfigure}\hfill
	\begin{subfigure}[t]{0.45\textwidth}
		\centering
		\includegraphics[width=.9\textwidth]{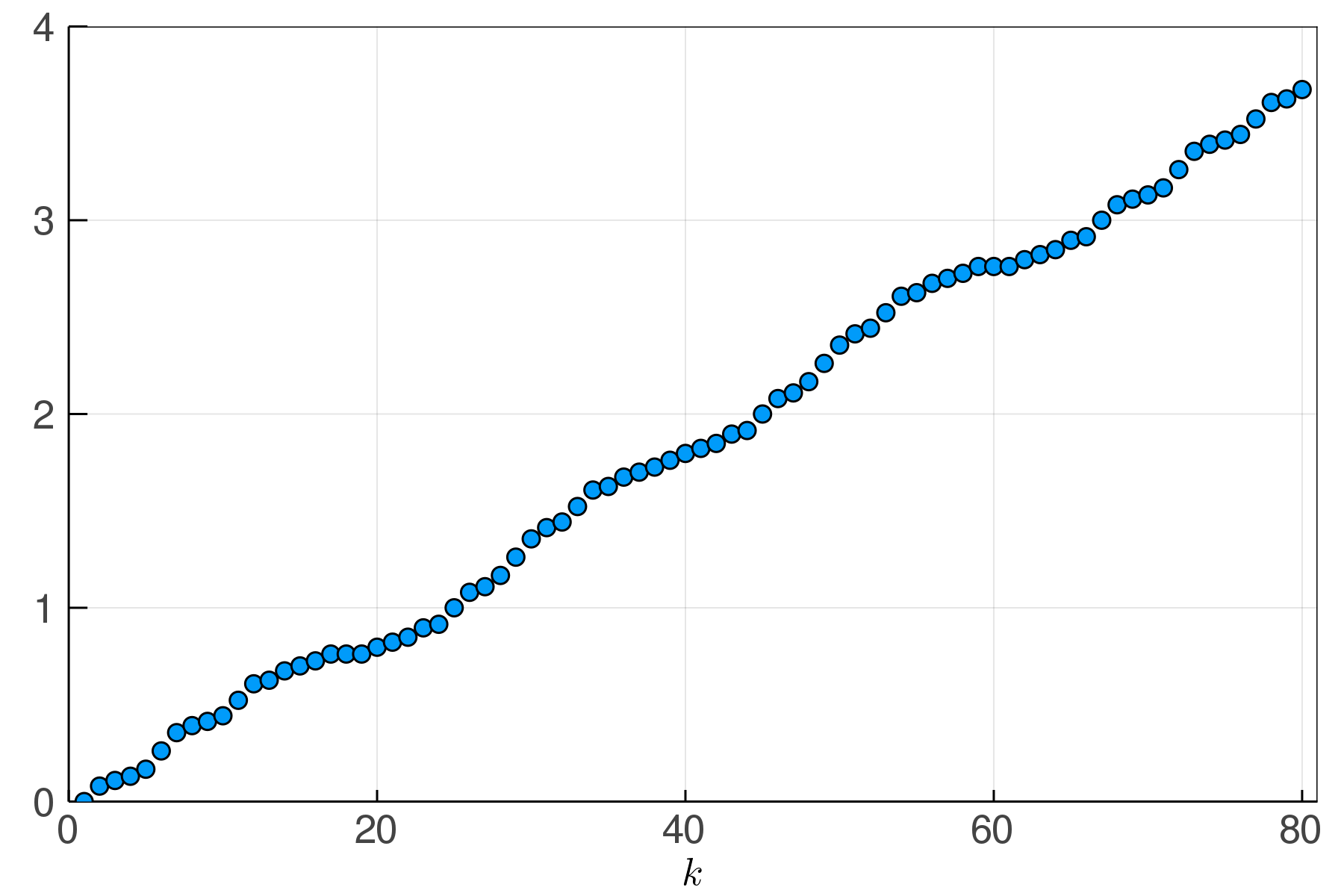}
		\caption{Lifted Coordinate.}
		\label{fig:coord-lifted}
	\end{subfigure}
	\caption{Coordinate change over time and lifted coordinate for the first $80$ time steps of the time series in Fig.~\ref{fig:trajectory_and_covering} and the coordinate in Fig.~\ref{fig:coord-good1}. One can see from Fig.~\ref{fig:coord-time} that the time series does approximately $3.5$ turns with respect to the coordinate. The lifted coordinate Fig.~\ref{fig:coord-lifted} captures this property, as it increases by $ 3.5 $.}
	\label{fig:normal-vs-lifted}
\end{figure}

\paragraph{2. Lifted Coordinates.}
Given a quantized time series $ \hat{x} : [1,T]\rightarrow Z$ and a circular coordinate $ \theta : Z \rightarrow S^1 $ we can form the composite $ \theta \circ\hat{x} : [1,T]\rightarrow S^1 $, which captures the change of the coordinate $\theta$ over time. Analogous to continuous maps, we lift this function to a function $ \hat{\theta} : [1,T]\rightarrow \R$ such that $ \pi_{S^1}\circ \hat{\theta} = \theta\circ \hat{x}$: We define the \emph{lifted coordinate} of $\theta$ and $\hat{x}$ via $ \hat{\theta}(1) = 0 $ and
\[
\hat{\theta}(t) = \hat{\theta}(t-1) + d_{S^1}(\theta\circ\hat{x}(t),\theta\circ\hat{x}(t-1))
\]
where $ d_{S^1}(x,y) $ denotes the signed geodesic distance from $y$ to $x$ on $S^1$. Figure \ref{fig:normal-vs-lifted} shows an example for a coordinate and its lift. 

\paragraph{3. Dynamically Relevant Coordinates}
We now identify coordinates which are relevant for dynamics. This is done in three steps: First, we define the correlation of coordinates which is large if coordinates describe the same features of the dynamics. Second, we search for a correlation-minimal basis of $  H^1(\cK,\Z)$. Finally, we discard all coordinates from this basis which do not describe any cycling dynamics.

We begin by defining the correlation number of two lifted coordinates $ \hat{\theta} $ and $ \hat{\eta} $ as
\begin{equation}
	c(\hat{\theta},\hat{\eta}) = \langle | \Delta \hat{\theta} |, |\Delta \hat{\eta}|\rangle
	\label{eq:correlation-two-coordinates}
\end{equation}
where the $i$-th entry of the vector $ \Delta \hat{\theta}$ is the forward finite difference $ \hat{\theta}_{i+1}-\hat{\theta}_i$ and $\langle\cdot,\cdot\rangle$ denotes the standard Euclidean scalar product.
For a basis $B $ of $ H^1(\cK,\Z) $ we define its \emph{correlation number} as 
\begin{equation}
	I\left( B\right) = \sum_{\substack{\alpha,\alpha'\in B\\\alpha \neq \alpha'}} c(\hat{\theta}_\alpha, \hat{\theta}_{\alpha'}).
	\label{eq:correlation-basis}
\end{equation}

These definitions can be motivated as follows: If two coordinates are to describe different features of the dynamics of a time series, they should change at disjoint periods of time, this is measured in \eqref{eq:correlation-two-coordinates}. For an optimal basis, we therefore minimize the overall correlation which is written out in \eqref{eq:correlation-basis}.

\begin{figure}[h]
	\centering
	\begin{subfigure}[t]{0.45\textwidth}
		\centering
		\includegraphics[width=.9\textwidth]{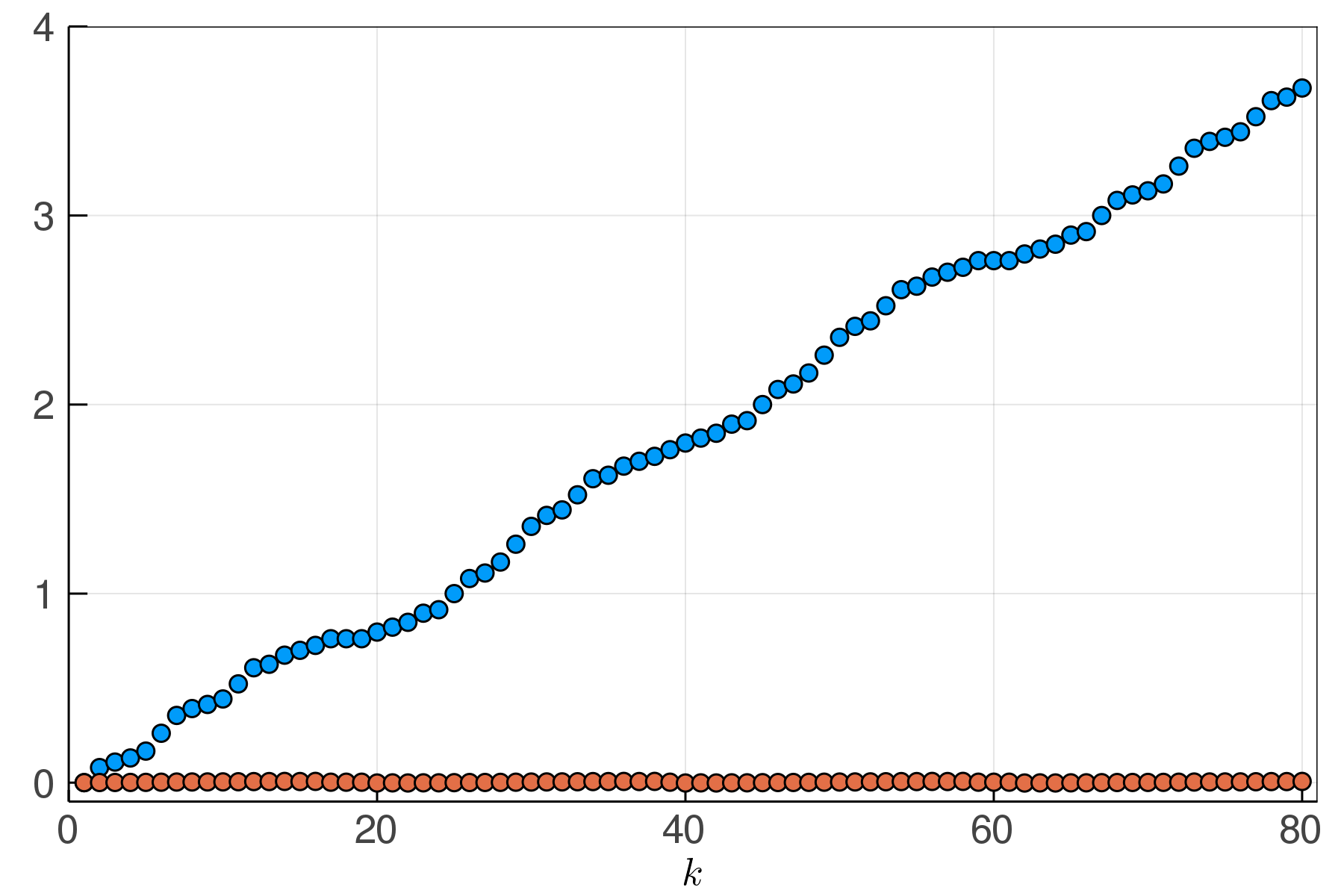}	
		\caption{Lift of the coordinates in Fig.~\ref{fig:coord-good}.}
		\label{fig:coord-lifted-good}
	\end{subfigure}\hfill
	\begin{subfigure}[t]{0.45\textwidth}
		\centering
		\includegraphics[width=.9\textwidth]{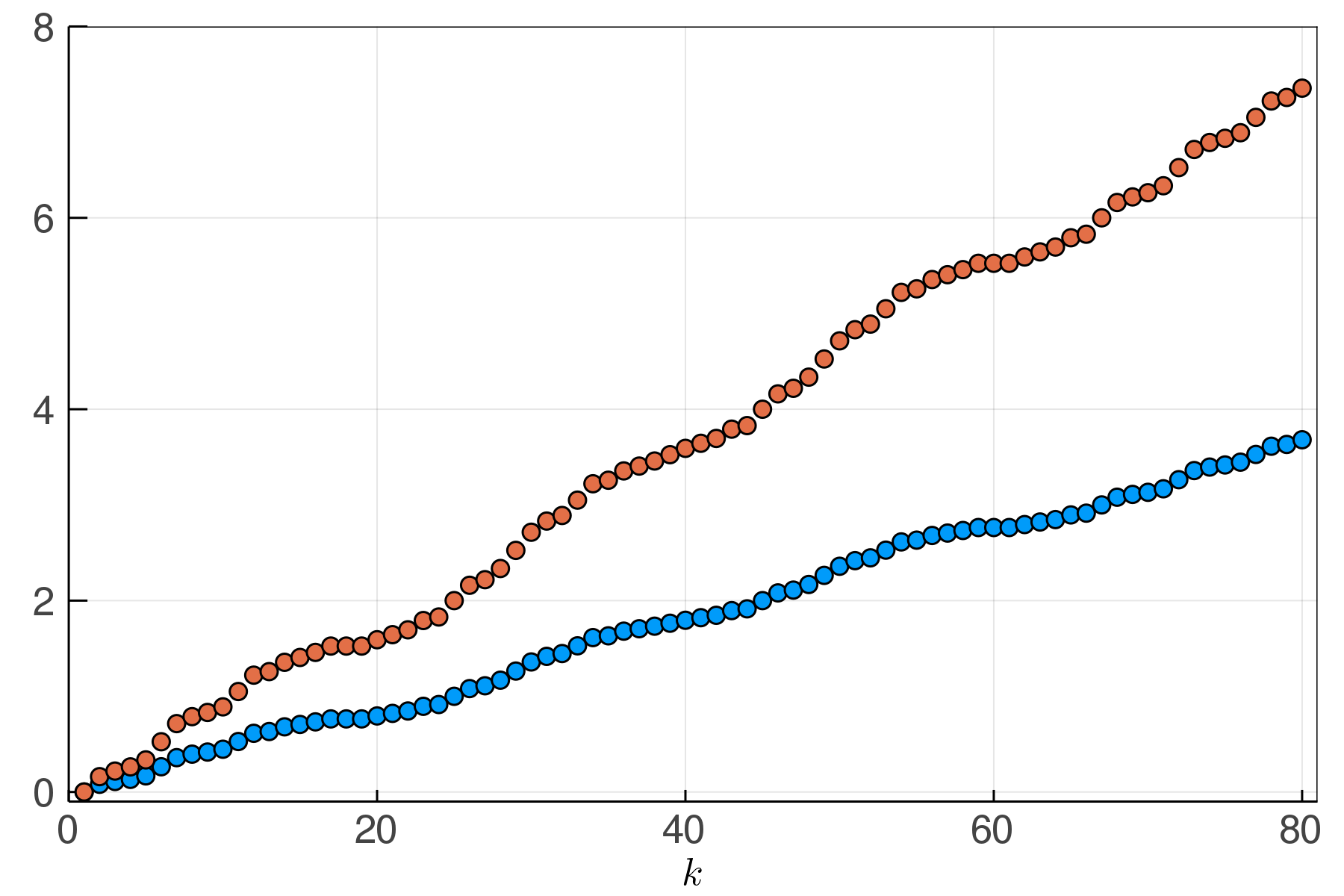}
		\caption{Lift of the coordinates in Fig.~\ref{fig:coord-bad}.}
		\label{fig:coord-lifted-bad}
	\end{subfigure}
	\caption{Lifted coordinates for two different basis.}
	\label{fig:lifted-good-vs-bad}
\end{figure}

As an example, we again consider the time series \ref{fig:quantization}. From Fig.~\ref{fig:coord-lifted} we know that for the first $80$ time steps, the series does $3.5$ turns around the right hole. 
Now consider Fig.~\ref{fig:lifted-good-vs-bad} where the lifted coordinates for the generating sets of figures Fig.~\ref{fig:coord-good} and Fig.~\ref{fig:coord-bad} are plotted for the first $80 $ time steps. The plots indicates that the lifted coordinates in Fig.~\ref{fig:coord-lifted-good} have a lower correlation number than the ones in Fig.~\ref{fig:coord-lifted-bad}. An explicit computation (for all $1000$ time steps) yields the values $0.158$ and $12.8 $, respectively, confirming that the preferred basis has lower correlation.

We now search for a basis with minimal correlation. Assuming $ \alpha_1,\dots,\alpha_n $ is any basis for the free group $ H^1(X,\Z) $, every basis can be written as $ A\alpha_1,\dots,A\alpha_n $ where $ A\in \operatorname{GL}_n(\Z) $. 
In order to find a correlation minimizing basis, we search $\operatorname{GL}_n(\Z) $, starting with the identity $ A := \operatorname{I}$ and recursively applying basis change operations (sums/swaps of rows/columns, multiplication of rows/columns with a unit) to $A$ up to a given depth. Of all these matrices we return the one with minimal correlation.

This approach works sufficiently well for simple examples. However, since we are only searching a finite subset of $ \operatorname{GL}_n(\Z) $ we have no guarantee of actually finding a minimizer (if one even exists). A better algorithm for finding a correlation minimal basis is a topic for future work.

At this point, the set $ \{\theta_\alpha\}_{\alpha\in B} $ contains $ \operatorname{dim}(B) $ many coordinates. Since the mere presence of a $1$-cycle does not imply the existence of cycling motion around the corresponding coordinate, we have to identify those coordinates in $ \{ \theta_\alpha\}_{\alpha\in B} $ which capture cycling behavior. For this, note that whenever the function $\hat \theta_\alpha$ is monotonic on some interval $[k,\ell]$, the trajectory is moving along the cyclic coordinate $\theta_\alpha$. We define the time series $ \hat{x} $ to be \emph{cycling along $\alpha$} if there is an interval $[k,\ell]$ where $ \hat{\theta}_\alpha $ satisfies a \emph{monotonicity criterion} and $|\hat \theta_\alpha(\ell)-\hat \theta_\alpha(k)| > 1$. The second condition ensures that the time series completes at least one full turn during the segment $ [k,\ell] $.
For a monotonicity criterion, a possible choice is 
\begin{center}
	$ |\theta_\alpha(t+1) - \theta_\alpha(t)|>\varepsilon$ for all $ t\in [k,\ell-1] $ and a fixed $ \varepsilon > 0 $.
\end{center}
$ \theta_\alpha $ is then said to be \emph{$\varepsilon$-increasing} along $ [k,\ell] $. In practice, we slightly relax this criterion and only require $ |\theta_\alpha(t+2) - \theta_\alpha(t)|>\varepsilon$ since sometimes adjacent cubes get assigned the exact same coordinate value. In this case, $\theta_\alpha$ is said to be \emph{almost $\varepsilon$-increasing} on $ [k,\ell] $. Note that the parameter $ \varepsilon$ has to be specified by the user; we typically do this by inspecting the lifted coordinates.

We define the subset $E\subset B $ of all \emph{dynamically relevant generators} of the basis as all $\alpha\in B$ for which the time series is cycling along $ \theta_\alpha $. The elements in $ B\setminus E $ will be called \emph{spurious} generators.

\subsection{Macro model}
We transfer the information on cycling motion back to the cubical covering:  A cube $\xi$ in the covering $\cX$ is \emph{$\alpha$-cycling} if the trajectory is cycling along $\alpha$ on some interval $[k,\ell]$ and there is $j\in [k,\ell]$ such that $\hat{x}_j\in\xi$. 
For $\xi\in\cX$, let $E(\xi)\subset E$ be the set of all dynamically relevant generators $\alpha$ for which $\xi$ is $\alpha$-cycling.

The cubical covering $\cX$ can now be decomposed into equivalence classes: Two cubes are equivalent if they are cycling with respect to the same set of non-spurious generators of $H^1(\cK;\Z)$:
\[
\xi \sim_1 \xi' \iff E(\xi)=E(\xi').
\]

We can furthermore distinguish cubes in which the trajectory ceases to be cycling. For this, assume the time series is $\alpha$-cycling along an interval $[k,\ell]$ which is maximal in the sense that the time series is not $\alpha$-cycling on any interval which contains $ [k,\ell] $. Now let $t$ be the first time step such that $ | \hat{\theta}_\alpha(\ell) - \hat{\theta}_\alpha(t) | <1 $. Then the cubes $Q^{-1}(\hat{x}_t),\dots,Q^{-1}(\hat{x}_\ell)$ are precisely those cubes which are hit during the last full turn with respect to $\alpha$ in $[k,\ell]$. We call such cubes $\alpha$-\emph{transient}. For a given cube $\xi$, we let $ E_t(\xi) $ denote the set of all generators which $\xi$ is transient for.

As a finer classification of cubes we can define 
\[
\xi \sim_2 \xi' \iff E(\xi)=E(\xi')\text{ and } E_t(\xi)= E_t(\xi').
\]

We now classify the cubes in $\cX$ according to either of these two equivalence relations and count transitions between the classes. That is, we build the quotient
\[
[\cX] := \cX/\sim \;= \{ [\xi_1]_\sim,\ldots,[\xi_q]_\sim \}
\]
as well as the transition matrix
\[
P(\sim) = (p_{ij}), \quad p_{ij} = \# \{ t \in [1,T-1] \mid \hat{x}_t \in [\xi_j]_\sim,\, \hat{x}_{t+1} \in [\xi_i]_\sim \}.
\]
We now call $(\cX/\sim_1, P(\sim_1))$ a macro model, and $(\cX/\sim_2, P(\sim_2))$ an extended macro model for the given time series. 

\begin{figure}[H]
	\centering
	\begin{subfigure}[b]{0.38\textwidth}
		\includegraphics[width=\textwidth]{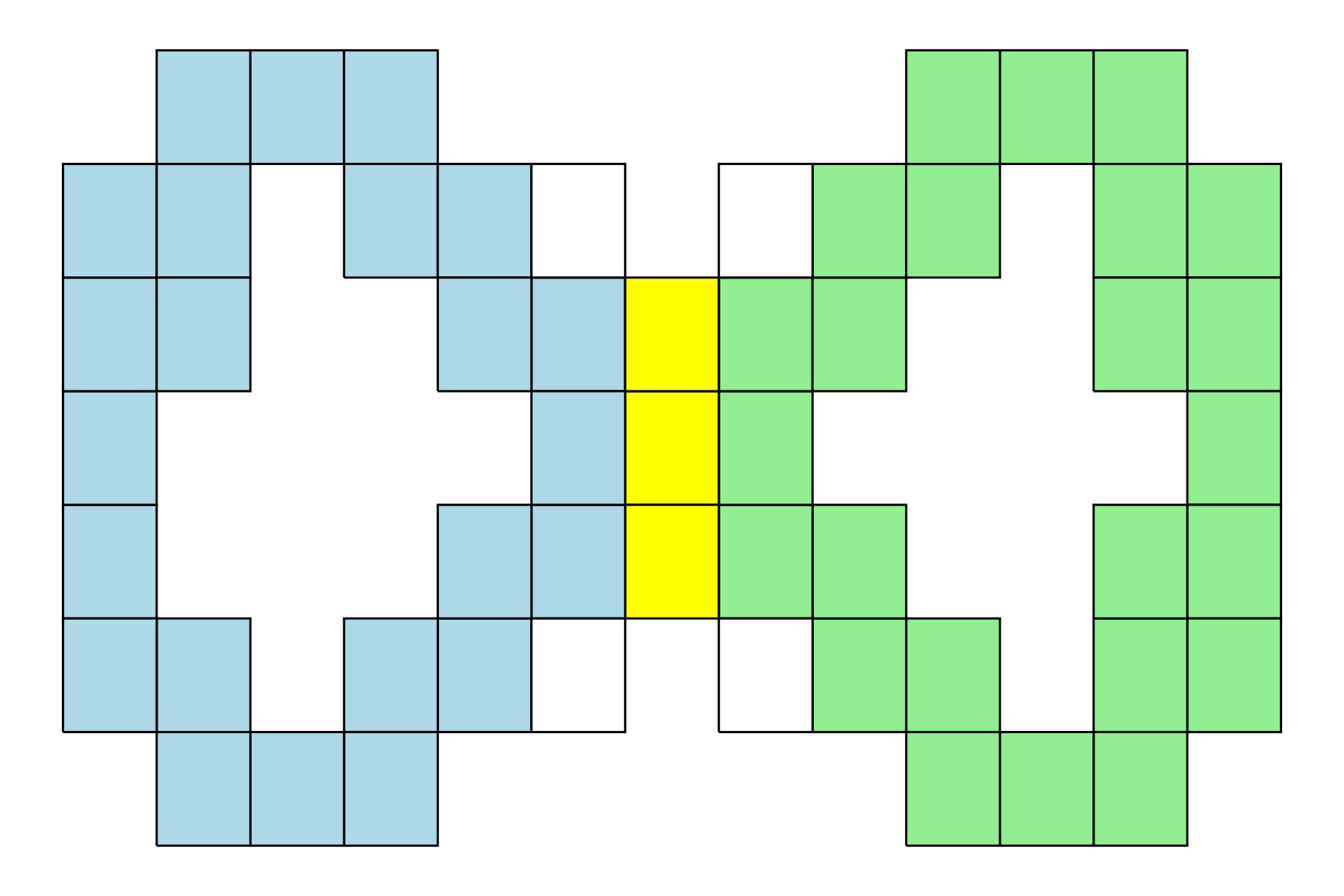}	
		\caption{decomposition}
	\end{subfigure}
	\begin{subfigure}[b]{0.58\textwidth}%
		\begin{tabular}{ r r r r l }
			\mmbox{julia_lightblue}  & \mmbox{julia_lightgreen} & \mmbox{julia_yellow} & \mmbox{julia_white} & \\ 
			1161 &  & 45 & 47 &\mmbox{julia_lightblue} $\alpha$-trn \\
			& 1264 & 55 & 49 &\mmbox{julia_lightgreen }$\beta$-trn \\
			44 & 56 &  35 &  &\mmbox{julia_yellow} $\alpha$-trn $\beta$-trn\\
			47 & 49 &  &  & \mmbox{julia_white} no cycling \\						
		\end{tabular}%
		\caption{transition matrix}
	\end{subfigure}
	\caption{Macro model for the double well system.}
	\label{fig:macro}
\end{figure}

\begin{figure}[H]
	\centering
	\begin{subfigure}[b]{0.38\textwidth}
		\includegraphics[width=\textwidth]{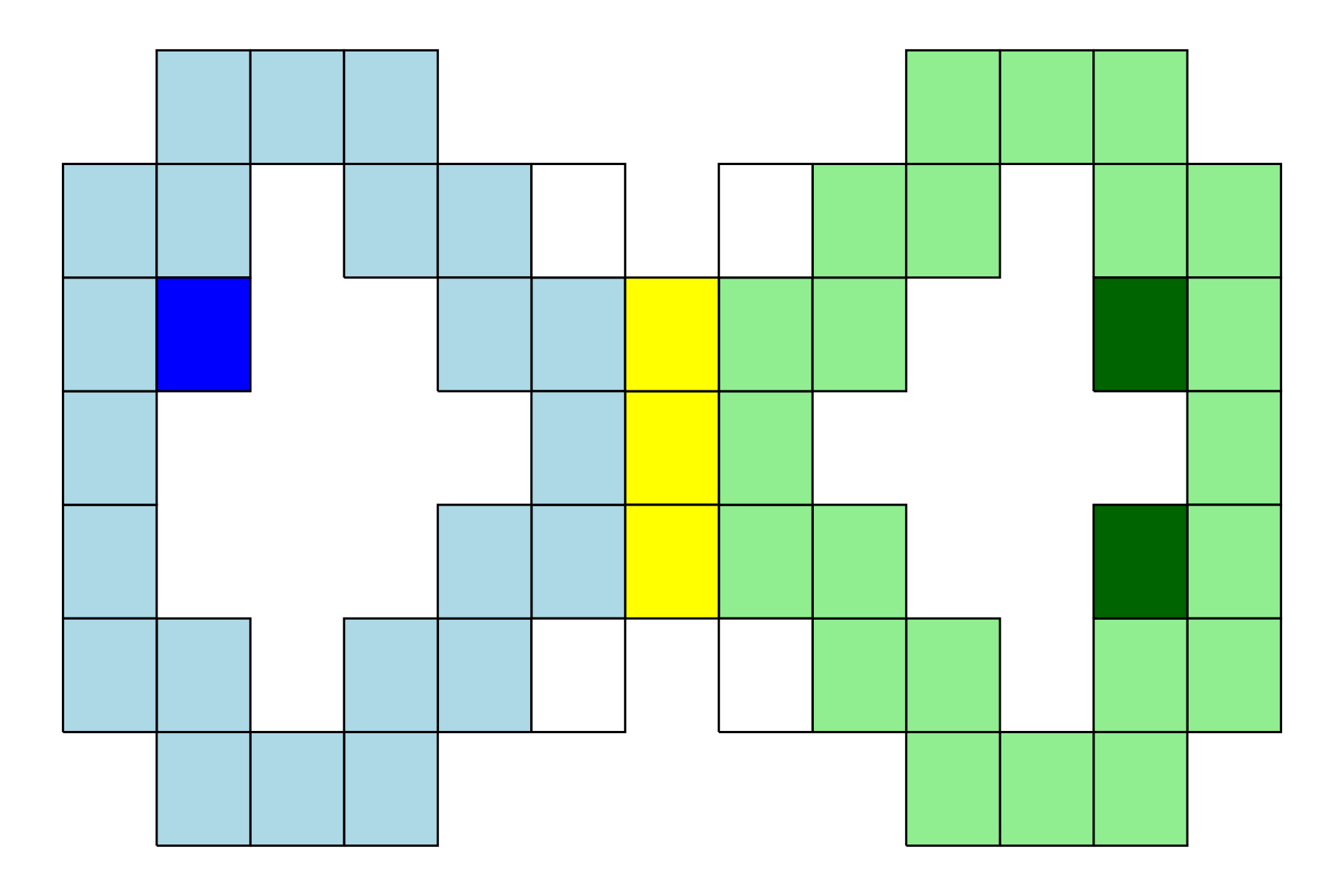}	
		\caption{decomposition}
	\end{subfigure}
	\begin{subfigure}[b]{0.58\textwidth}%
		\begin{tabular}{ r r r r r r l }
			\mmbox{julia_darkblue} & \mmbox{julia_lightblue} & \mmbox{julia_darkgreen} & \mmbox{julia_lightgreen} & \mmbox{julia_yellow} & \mmbox{julia_white}\\
			   & 4 &   &   &   &   & \mmbox{julia_darkblue} $\alpha$-cyc \\
			 4 & 1256 &   &   & 55 & 49 & \mmbox{julia_lightblue} $\alpha$-trn \\
			   &   &   & 4 &   &   & \mmbox{julia_darkgreen} $\beta$-cyc \\
			   &   & 4 & 1153 & 45 & 47 & \mmbox{julia_lightgreen} $\beta$-trn \\
			   & 56 &   & 44 & 35 &   & \mmbox{julia_yellow} $\alpha$-trn $\beta$-trn \\
			   & 49 &   & 47 &   &   & \mmbox{julia_white} no cycling \\ 
		\end{tabular}%
		\caption{transition matrix}
	\end{subfigure}
	\caption{Extended macro model for the double well system.}
	\label{fig:macro-with-transient}
\end{figure}

Figs.~\ref{fig:macro} and \ref{fig:macro-with-transient} show both macro models for the double well example.


We remark that by rescaling the columns of $P$ one obtains a Markov matrix with entries which could be interpreted as transition probabilities between different regions. We deliberately choose not to define transition matrices this way, since transitions generally do not satisfy a Markov property. 



\section{Experiments}

The following results are obtained using our implementation of the pipeline in Section~\ref{sec:construction} in the programming language \href{https://julialang.org}{Julia} \cite{BeEdKaVi17}. In particular, we use the algorithm in \cite{HaMiMrNa14} for computing $H^1$ with integer coefficients.

\subsection{Perturbed double well}
\label{sec:dw-example}
The time series in Fig.~\ref{fig:trajectory_and_covering}, which was used to illustrate the constructions in Section \ref{sec:construction}, was obtained by integrating a stochastically perturbed version of the double well Hamiltonian system
\begin{equation}
	\label{eq:SDE}
	dx = f(x) dt + \sigma dB,
\end{equation}
with $x=(q,p)$, $f(x) = (p, q-q^3)$, $\sigma = (0,0.025)$ and $B$ denoting Brownian motion. We integrate \eqref{eq:SDE} from the initial value $x=(1,0.7)$ by the SRIW1 scheme \cite{Ro10} with step size 0.01.

The macro models in Fig.~\ref{fig:macro} and Fig.~\ref{fig:macro-with-transient} were generated using the coordinates in Fig.~\ref{fig:coord-good} and the monotonicity criterion ''$ \theta $ is almost $ \varepsilon $-increasing'' with $\varepsilon= 0.04$. We note that Fig.~\ref{fig:coord-lifted-good} provides a hint for choosing $\varepsilon$ since 3.5 turns in 80 steps average to an increase of approximately $ 0.04 $ per step. 
These models are like we would expect for such simple dynamics. The yellow boxes capture the location in phase space where direct transitions between the loops is possible and the green and blue boxes capture the location where the trajectory cycles around the natural holes.



\subsection{The Lorenz system}

For this example we generated a time series by integrating the Lorenz system with the classical parameters $\sigma =10,\, \beta = \frac{8}{3} $ and $\rho=28$ with time step size $0.1$ for $1$ million time steps using the classical fourth order Runge Kutta method. As starting value, we choose $(0,10,0)$, but we discard the first 6000 time steps since they 'close up' the left holes of the complex. This highlights one shortcoming of our current technique which will be addressed in future work.

\begin{figure}[H]
	\centering
	\begin{subfigure}[b]{0.23\textwidth}
		\centering
		\includegraphics[width=\textwidth]{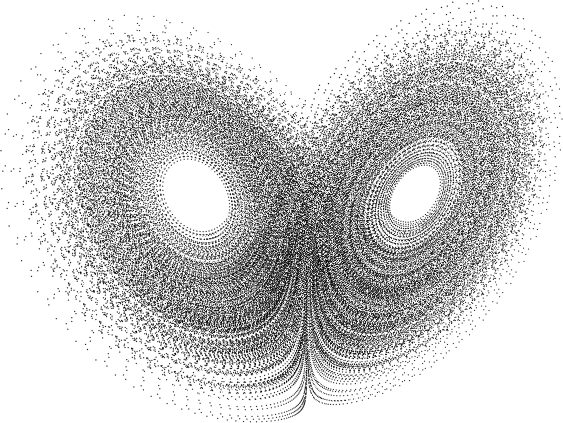}
		\caption{Sampled trajectory}
		\label{fig:lorenz-trajectory}
	\end{subfigure}\hfill
	\begin{subfigure}[b]{0.23\textwidth}   
		\centering
		\includegraphics[width=\textwidth]{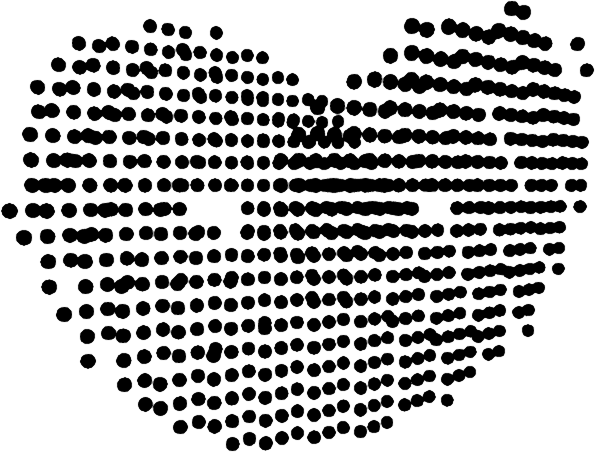}
		\caption{Box states}
		\label{fig:lorenz-quantization}
	\end{subfigure}\hfill
	\begin{subfigure}[b]{0.23\textwidth}  
		\centering 
		\includegraphics[width=\textwidth]{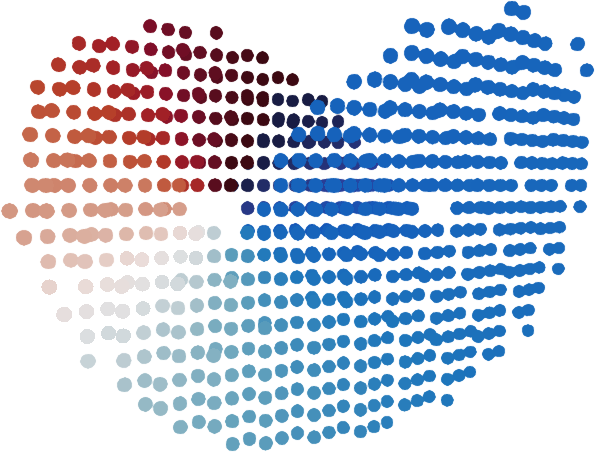}
		\caption{First coordinate}
		\label{fig:lorenz-coord1}
	\end{subfigure}\hfill
	\begin{subfigure}[b]{0.23\textwidth}   
		\centering
		\includegraphics[width=\textwidth]{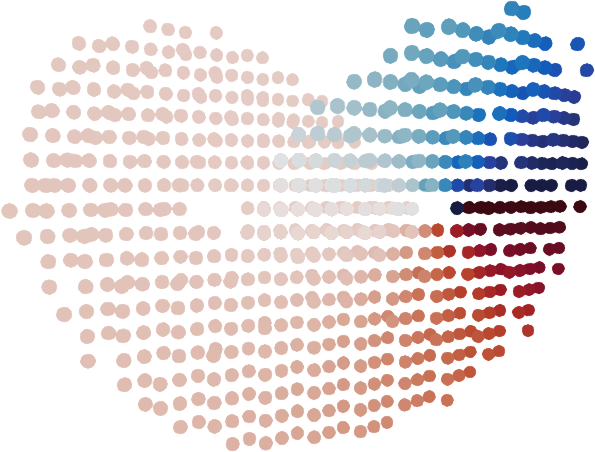}
		\caption{Second coordinate}
		\label{fig:lorenz-coord2}
	\end{subfigure}
	\caption{Illustration of the pipeline for a trajectory on the Lorenz attractor. }
	\label{fig:lorenz-example}
\end{figure}

We choose the quantization radius $ r=2.5 $ and construct a cubical covering with 652 cubes and a quantized time series with 42698 points. Since the first cohomology of the resulting Vietoris-Rips complex is two-dimensional, we compute two circular coordinates, which are both seen to be dynamically relevant with the monotonicity criterion ''$\theta $ is almost 0.007-increasing''. Plots for the time series and the intermediate steps can be found in Fig.~\ref{fig:lorenz-example}. 
\begin{figure}[H]
	\centering
	\begin{subfigure}[b]{0.35\textwidth}
		\includegraphics[width=\textwidth]{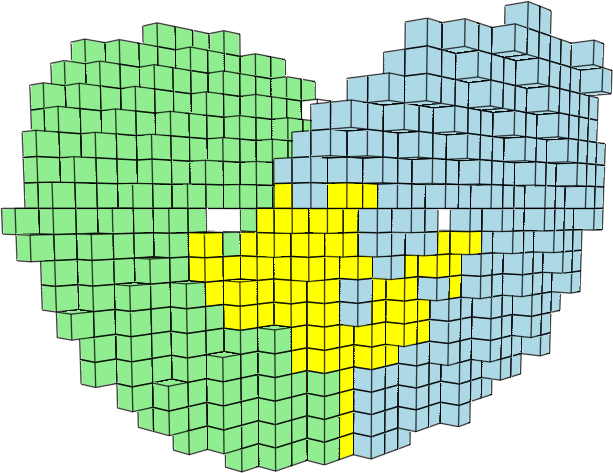}	
		\caption{Decomposition of box covering.}
		\label{fig:lorenz-decomposition}
	\end{subfigure}%
	\begin{subfigure}[b]{0.59\textwidth}%
		\begin{tabular}{ r r r r l }
			\mmbox{julia_lightblue}  & \mmbox{julia_lightgreen} & \mmbox{julia_yellow} & \mmbox{julia_white} & \\ 
			15925 &  & 818 &  &\mmbox{julia_lightblue} $\alpha$-cycling \\
			& 15011 & 876 & 2 &\mmbox{julia_lightgreen } $\beta$-cycling \\
			818 & 875 & 8370 &  &\mmbox{julia_yellow} $\alpha$-cycling \& $\beta$-cycling\\
			 & 2 &  &  & \mmbox{julia_white} no cycling \\						
		\end{tabular}%
		\caption{Transition matrix}
	\end{subfigure}
	\caption{Macro model for the Lorenz system. The cubical covering with $652$ boxes is decomposed into 277 $\alpha$-cycling, 256 $\beta$-cycling and 108 $\alpha$-$\beta$-cycling cubes as well as 2 cubes which do not contain any cycling motion.}
	\label{fig:lorenz-macro}
\end{figure}
The macro model shown in Fig.~\ref{fig:lorenz-macro}  nicely captures many important aspects of the dynamics on the Lorenz attractor. We learn that there are (at least) two different types of cycling motion, that each of these occurs in a distinct region in phase space (the blue and green regions), and that these regions intersect (yellow region). We furthermore see that cycling dynamics are present almost everywhere in the box decomposition as there are only 2 non-cycling boxes. 

\begin{figure}[H]
	\centering
	\begin{subfigure}[b]{0.35\textwidth}
		\includegraphics[width=\textwidth]{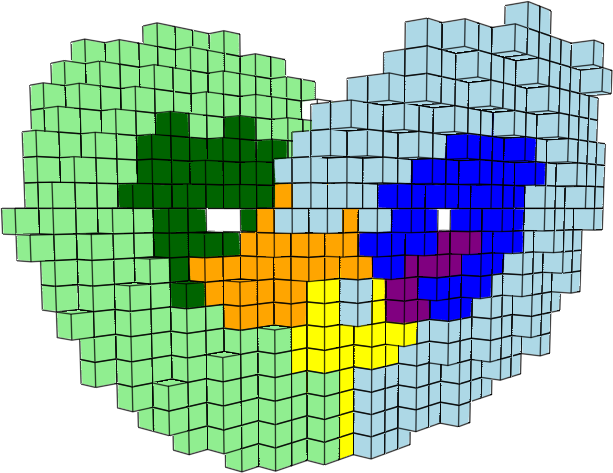}	
		\caption{decomposition}
		\label{fig:lorenz-decomposition-with-transient}
	\end{subfigure}\\
	\begin{subfigure}[b]{0.59\textwidth}%
		\begin{tabular}{ r r r r r r r r l }
			\mmbox{julia_darkblue} & \mmbox{julia_lightblue} & \mmbox{julia_darkgreen} & \mmbox{julia_lightgreen} & \mmbox{julia_yellow} & \mmbox{julia_orange} & \mmbox{julia_purple} & \mmbox{julia_white} & \\ 
			 4065 & 1000 &   &   &   & 11 & 263 &   & \mmbox{julia_darkblue} $\alpha$-cyc\\
			 968 & 9892 &   &   & 80 & 464 &   &   & \mmbox{julia_lightblue} $\alpha$-trns\\
			   &   & 2953 & 789 &   & 222 & 8 &   & \mmbox{julia_darkgreen} $\beta$-cyc\\
			   &   & 701 & 10568 & 92 & 82 & 472 & 2 & \mmbox{julia_lightgreen} $\beta$-trns\\
			   & 366 &   & 401 & 3133 & 160 & 301 &   & \mmbox{julia_yellow} $\alpha$-trns $\beta$-trns\\
			 14 & 106 & 306 & 66 & 452 & 2098 & 33 &   & \mmbox{julia_orange} $\alpha$-trns $\beta$-cyc \\
			 292 & 40 & 12 & 90 & 604 & 39 & 1550 &   & \mmbox{julia_purple} $\alpha$-cyc $\beta$-trns\\
			   &   &   & 2 &   &   &   &   & \mmbox{julia_white} no cycling\\\						
		\end{tabular}
		\caption{transition matrix}
	\end{subfigure}
	\caption{Extended macro model for the Lorenz system.}
	\label{fig:lorenz-macro-with-transient}
\end{figure}
In the extended macro model (Fig.~\ref{fig:lorenz-macro-with-transient}), we see that the cycling regions are subdivided into a cycling set near the center of the wings and a transient set near the outside of the wings. This indicates that all cycling dynamics in the inside of the wings eventually moves to the outer regions. The extended model furthermore identifies the regions where direct transitions between cycling dynamics can occur. The purple and yellow regions in Fig.~\ref{fig:lorenz-decomposition-with-transient} are the only regions where a direct transition from $\alpha$- to $\beta$-cycling dynamics is possible and the orange and yellow regions are the only places that can contain the reverse transition.




%

\section{Discussion}

The techniques described in this paper appear to be a promising novel approach to identifying from time series data regions of phase space in which oscillations occur and locations at which transitions between these oscillations occur. However, a number of distinct questions need to be answered to obtain confidence in applying this technique to complicated higher-dimensional systems, where the results cannot be inspected and modified by visualization. We briefly address those in the following paragraphs.

\paragraph{Construction of the complex.} The computation of circle-valued coordinates from data requires the construction of a geometric complex. The approach chosen in this article accomplishes this by constructing a Vietoris--Rips complex from a suitably quantized version of the given time series. In particular, we rely on finding a quantization radius $r$ which is small enough to contain those holes which give rise to dynamically relevant coordinates and large enough to connect the data in a meaningful way. In general, such a radius need not exist. This even happens in the Lorenz system with a trajectory that starts very close to the center of one of the wings.

\paragraph{Finding optimal coordinates.} When searching for dynamically relevant coordinates, we encounter the problem of finding a correlation minimizing basis. This poses the natural question of existence and uniqueness of such a basis. In addition, an algorithm is needed to compute this basis or a suitable approximation. Furthermore, since sparse data can lead to a large dimension of $H^1(\cdot,\Z)$, an efficient algorithm for this computation would be desirable. We hope to be able to address this using techniques inspired by those for the computation of (persistent)  cohomology with coefficients in a finite field \cite{Ba19}.

\paragraph{Identifying cycling motion.} In this contribution, cycling motion is identified by analyzing the monotonicity behavior of circle valued coordinates. While this leads to satisfactory results in the presented examples, we have no general reliable procedure of identifying recurrence. For example, a cycling time series with a bit of back-and-forth moving in every full turn would be difficult to identify using the presented methods.






\end{document}